\documentclass[12pt,leqno]{article}
\usepackage{a4wide}
\usepackage{amssymb,amsmath}
\allowdisplaybreaks[1]

\newtheorem{theorem}{Theorem}%[section]
\newtheorem{corollary}[theorem]{Corollary}
\newtheorem{definition}[theorem]{Definition}

\newtheorem{lemma}[theorem]{Lemma}

\newenvironment{proof}{\begin{trivlist}\item[]{\it
Proof.}}{\hfill$\square$\end{trivlist}}
\newcommand{\Q}{\bigskip\par\noindent}

\def\widebar{\overline}

\def\gkdim{{\rm GKdim}}
\def\span{{\rm span}}
\def\mn{{\mathbb N}}

\begin{document}
\author{{\bf T H Lenagan}\\
{\small\em School of Mathematics, University of Edinburgh}\\
{\small\em James Clerk Maxwell Building, King's Buildings}\\
{\small\em Mayfield Road, Edinburgh EH9 3JZ, Scotland}\\
{\small\em tom@maths.ed.ac.uk}\\
{\bf Agata Smoktunowicz}\\
{\small\em Institute of Mathematics, Polish Academy of Sciences}\\
{\small\em Sniadeckich 8, 00-956  Warsaw 10, Poland}\\ {\small\em
E-mail: agatasm@impan.gov.pl}}

\title{An infinite dimensional affine nil algebra with finite
Gelfand-Kirillov dimension\thanks{Part of this work was done while the
second author was visiting the University of Edinburgh, with support
from the Edinburgh Mathematical Society. The first author acknowledges
support by Leverhulme Grant F/00158/X.}} 

\date{ } 

\maketitle

\begin{abstract} We construct a nil algebra  over a countable field
 which has finite but non-zero Gelfand-Kirillov dimension.
 \end{abstract}


\section*{Introduction} 

In their famous paper \cite{gs}, Golod and Shafarevich gave a
construction which yields, over arbitrary fields, infinite
dimensional, finitely generated nil algebras that are not nilpotent
algebras. In fact, the algebras arising from the Golod-Shafarevich
construction have exponential growth.  In view
of this, Lance Small has asked whether there is a finitely generated,
infinite dimensional nil algebra with finite Gelfand-Kirillov
dimension, which is not nilpotent, \cite{small}.  Small's question is
answered in this paper: we construct, over an arbitrary countable
field, a finitely generated, non-nilpotent, nil algebra with
Gelfand-Kirillov dimension not exceeding $20$.

In what follows $K$ is a countable field and $A$ is the free $K$-algebra 
in three non-commuting indeterminates $x, y$ and $z$.  The set of
monomials in $x$, $y$, $z$ is denoted by $M$ and $M(n)$ denotes the
set of monomials of degree $n$, for each $n\geq 0$.  Thus,
$M(0)=\{1\}$ and for $n\geq1$ the elements in $M(n)$ are of the form
$x_{1}\ldots x_{n}$, where $x_{i}\in \{x, y, z\}$. The $K$-subspace of
$A$ spanned by $M(n)$ will be denoted by $H(n)$ and elements of $H(n)$
will be called {\em homogenous polynomials of degree $n$}.  Every
polynomial $f\in A$ such that $\deg(f)=d$ can be uniquely presented in
the form $f=f_{0}+f_{1}+\ldots +f_{d}$, where $f_{i}\in H(i)$.
The elements $f_{i}$ are the {\em  homogeneous components} of $f$ and
$\deg(f)$ denotes the degree of the polynomial $f$.  A right ideal $I$
of $A$ is {\em homogeneous} if for every $f\in I$ all homogeneous
components of $f$ are in $I$.  Let $V$ be a linear space over $K$,
then $\dim_{K}V$ denotes the dimension of $V$ over $K$.  The
Gelfand-Kirillov dimension of an algebra $R$ is denoted
by $\gkdim(R)$. For elementary properties of Gelfand-Kirillov dimension
we refer to \cite{kl}.

%%%%%%%%%%%%%%%%%%%%%%%%%%%%%%%%%%%%%%%%

\section{Enumerating elements} \label{one} 

Let $\widebar{A}$ be the subalgebra of $A$ consisting of polynomials
with constant term equal to zero. As usual, $\mn$ denote the set of natural
numbers. 

The aim is to present an algebra with the desired properties as
$\widebar{A}/E$ for a suitable ideal $E$. The ideas develop from ideas
in earlier papers by the second author, but we have to define several
sets of subspaces carefully in order to be able to control the growth
of the algebra we have in mind.

We start with two results derived from similar results in \cite{smok}. 

%%%%%%%%%%%%%%%%%%%%%%%%%%%%%%%%%%%%%%%%%%

\begin{lemma}\label{enumerate} 
 Let $K$ be a countable field, and let $\widebar{A}$ be as above.
 Then there exist a set $Z\subseteq \mn$, with all $i\in Z$ being
 greater than or equal to $5$,
such that elements of
 $\widebar{A}$ can be enumerated as $f_{i}$ for $i\in Z$ (that is, 
 $\widebar {A}=\{f_{i}\}_{i\in Z}$) and such that 
 $2^{2^{i}}>3^{6t_{i}}$ for each $i\in Z$, 
where $t_{i}$ is the degree of $f_{i}$.
\end{lemma}

\begin{proof}  The field $K$ is countable and the algebra $A$ is
finitely generated over $K$, so the elements of $\widebar {A}$ can be
enumerated: say $\widebar {A} = \{g_1, g_2, \dots \}$. We now define
an increasing function $\theta: \mn \longrightarrow \mn$ as
follows. 
Set $\theta(1):= \min \{ i \in\mn\mid i>4, \; 2^{2^{i}}>3^{6\deg(g_1)}\}$. 
Suppose that we have defined $\theta:\{1, \dots, n\}
\longrightarrow \mn$ such that $2^{2^{\theta(i)}} > 3^{6\deg(g_i)}$,
for each $i=1, \dots, n$. Then set $\theta(n+1) := \min\{s\mid s>
\theta(1), \dots, \theta(n) {\rm ~and~} 2^{2^s} >
3^{6\deg(g_{n+1})}\}$. If we now rename the elements of $\widebar{A}$ by
setting $f_{\theta(s)} := g_s$ then we have a listing of the elements
of $\widebar{A}$ with the required properties.
\end{proof}

%%%%%%%%%%%%%%%%%%%%%%%%%%%%%%%%%%%%%

Given a subset $S\subseteq H(n)$, for some $n$, let $B_{n}(S)$ denote
the right ideal of $A$ generated by the set
$\bigcup_{k=0}^{\infty}M(nk)S$; that is,
$$B_{n}(S)=\sum_{k=0}^{\infty}M(nk)SA.$$

%%%%%%%%%%%%%%%%%%%%%%%%%%%%%%%%%%%%%%

\begin{theorem}\label{existFi}
 Let $Z$ , $\{f_{i}\}_{i\in Z}$ be as in Lemma 1. Let $i\in Z$, and
 let $I$ be the two-sided ideal generated by ${f_{i}}^{10w_{i}}$ where
 $w_{i}=4.2^{2^{i}}$. Then there is a linear $K$-space $F_{i}\subseteq
 H(2^{2^{i}})$ such that $I\subseteq B_{w_{i}}(F_{i})$ and
 $\dim_{K}(F_{i})<2^{2^{i+1}}-2$.
\end{theorem}

\begin{proof} 
Apply \cite[Theorem 2]{smok} with  $f=f_{i}$, $r=2^{2^{i}}$,
$w=w_{i}=4.2^{2^{i}}$, and put $F_{i}=\span_{K}F$.  Note that these choices of 
$f,r,w$
satisfy the hypotheses of \cite[Theorem 2]{smok}, by Lemma~\ref{enumerate}. 
\end{proof} 

%%%%%%%%%%%%%%%%%%%%%%%%%%%%%%%%%%%%%%%%

\section{Definition of $U(2^n)$ and $V(2^n)$}

Set $S:=\{[2^{i}-i-1, 2^{i}-1] \mid i=5, 6, 7,  \ldots \}$.

\begin{theorem} \label{7props} 
Let $Z$, $F_{i}$  be as in Theorem~\ref{existFi}.
 Then there are $K$-linear subspaces 
 $U(2^{n})$ and $V(2^{n})$ of $H(2^{n})$ such that for all $n>0$ we have:
 \begin{enumerate}
\item $\dim_{K}V({2^{n}})=2$ if $n\notin S$.
\item $\dim_{K}V(2^{2^{i}-i-1+j})=2^{2^{j}}$, for all $1<i$ and
all $0\leq j\leq i$.
\item $V(2^{n})$ is generated by monomials.
\item $F_{i}\subseteq U(2^{2^{i}})$ for every $i\in Z$.
\item $V(2^{n})+U(2^{n})=H(2^{n})$ and $V(2^{n})\cap U(2^{n})=0$.
\item $H(2^{n})U(2^{n})+U(2^{n})H(2^{n})\subseteq U(2^{n+1})$ .
\item $V(2^{n+1})\subseteq V(2^{n})V(2^{n})$.
 \end{enumerate}
\end{theorem}

\begin{proof} We construct the sets $U(2^n)$ and $V(2^n)$ inductively. 
Set $V(2^{0}):=V(1)=
Kx+Ky$ and $U(2^{0})=U(1):=Kz$.  Assume that we have 
defined
$V(2^{m})$ and  $U(2^{m})$ for $m\leq n$ in such a way that 
conditions 1-5 hold for all $m\leq n$ and conditions 6 and 7 hold for all 
$m<n$.  
Then we define $V(2^{n+1})$ and  $U(2^{n+1})$ in the following way.
Observe first that since $U(n)\cap V(n)=0$ then
  $$
\left\{
U(n)U(n)+U(n)V(n)+V(n)U(n)\right\}
\cap 
\left\{
V(n)V(n)
\right\} 
=0.
$$ 
Our next step is to make the following observation. If $\widebar {V},
P\subseteq V(n)V(n)$ and $\widebar {V}\cap P=0$ then
$$
\left\{
U(n)U(n)+U(n)V(n)+V(n)U(n)+\widebar {V}
\right\} 
\cap P=0.$$
For, suppose that  $c=c_{1}+c_{2} \in P$ with $c_{1}\in U(n)U(n)
 +U(n)V(n)+V(n)U(n)$ and $c_{2}\in \widebar {V}$.  
We claim that $c=0$.  Notice 
 that $c\in P$ and $c_{2}\in \widebar {V}$ implies that 
  $c_1 = c-c_{2}\in P+\widebar {V}\subseteq V(n)V(n)$.
  On the other hand,
   $c-c_{2}=c_{1}\in U(n)U(n)
 +U(n)V(n)+V(n)U(n)$. By the above observation, we get 
 $c_1=0$ so that $c=c_2\in \widebar{V}$. 
However, $c\in P$; so that $c\in P\cap \widebar{V}=0$, as required.

  Now we will define $V(2^{n+1})$,  $U(2^{n+1})$ inductively, in
  the following way.
  Consider the  three cases
  \begin{itemize}
  \item[1.] $n\in S$ and $n+1\in S$.
  \item[2.] $n\notin S$.
  \item[3.] $n\in S$ and $n+1\notin S$.
  \end{itemize}

\Q{\bf Case 1.} Suppose that $n\in S$ and $n+1\in S$.  Then we define
$V(2^{n+1}):=V(2^{n})V(2^{n})$; so condition 7 certainly holds.  Notice
that $V(2^{n+1})$ is spanned by monomials, since $V(2^{n})$ is spanned
by monomials; so condition 3 holds.  Moreover
$\dim_{K}V(2^{n})=(\dim_{K}V(2^{n}))^{2}$. Since $n, n+1\in S$, it
follows that $n=2^{i}-i-1+j$ for some $i$ and some $0\leq j<i$. By the
inductive hypothesis, $\dim_{K}V(2^{n})=2^{2^{j}}$.  Now
$\dim_{K}V(2^{n+1})=(2^{2^{j}})^{2}=2^{2^{j+1}}$, as required for
condition 2 (condition 1 does not apply in this case).  Set
$U(2^{n+1}):=U(2^{n})H(2^{n})+H(2^{n})U(2^{n})$; so condition 6
certainly holds. It is now easy to check that condition 5
holds. Finally, observe that since $n+1\in S$, we have $2^{n+1}\neq
2^{2^{i}}$ for every $i$, hence condition 4 is empty in this case and
so holds trivially.

\Q {\bf Case 2.} Suppose that  
$n\notin S$. Then $\dim_{K}V(2^{n})=2$, by the
inductive hypothesis.  Let $m_{1}, m_{2}$ be distinct monomials from
$V(2^{n})V(2^{n})$. Set $V(2^{n+1}):=Km_{1}+Km_{2}$.
Then $\dim_{K}V(2^n)=2$, as required.  Let $\widebar
{V}\subseteq V(2^{n})V(2^{n})$ be such that $\widebar {V}\cap
V(2^{n+1})=0$ and $\widebar {V}+V(2^{n+1})=V(2^{n})V(2^{n})$. Set 
$U(2^{n+1}):=U(2^{n})V(2^{n})+V(2^{n})U(2^{n})+U(2^{n})U(2^{n})+\widebar
{V}$.  We see that $U(2^{n+1})\cap V(2^{n+1})=0$ and
$U(2^{n+1})+V(2^{n+1})=H(2^{n+1})$.  Observe that since $n\notin S$,
we get $2^{n+1}\neq 2^{2^{i}}$ for every $i$, and again 4 holds trivially.

\Q 
{\bf Case 3.} Suppose that 
$n\in S$ while $n+1\notin S$.  Then $n=2^{i}-1$ for some
$i>1$.  By the inductive hypothesis 
$\dim_{K}V(2^{n})=\dim_{K}V(2^{2^{i}-1})=\dim_{K}V(2^{2^{i}-i-1+i})=2^{2^{i}}$.
Now $\dim_{K}V(2^{n})V(2^{n})=2^{2^{i+1}}$. 

Assume first that $i\in
Z$.  We know that $F_{i}$ has a basis $\{f_{1}, \ldots ,f_{s}\}$ for some
$f_{1},\ldots , f_{s}\in H(2^{2^{i}})$ and $s<2^{2^{i+1}}-2$.  Write
each $f_{j}$ as $f_{j}=\widebar {f}_{j}+g_{j}$ where $\widebar
{f}_{j}\in V(2^{2^{i}})V(2^{2^{i}})$ and $g_{j}\in
V(2^{2^{i}})U(2^{2^{i}})
+U(2^{2^{i}})U(2^{2^{i}})+U(2^{2^{i}})V(2^{2^{i}})$.  Since 
$V(2^{2^{i}})\cap U(2^{2^{i}})=0$ this decomposition is unique.  Let
$P$ be a $K$-linear subspace of $V(2^{2^{i}})V(2^{2^{i}})$ such that
$\widebar{f}_{j}\in P$ for all $1\leq j\leq s$ and
$\dim_{K}P=2^{2^{i+1}}-2$.

Since $V(2^{2^{i}})V(2^{2^{i}})$ is spanned by monomials and
 $\dim_{K}(V(2^{2^{i}})V(2^{2^{i}}))=2^{2^{i+1}}$ while
 $\dim_{K}P=2^{2^{i+1}}-2$, then there are monomials $m_{1}, m_{2}\in
 V(2^{2^{i}})V(2^{2^{i}}) $ such that
 $Km_{1}+Km_{2}+P=V(2^{2^{i}})V(2^{2^{i}})$ and $P\cap
 (Km_{1}+Km_{2})=0$.  Now set $V(2^{n+1}):=Km_{1}+Km_{2}$ and
 $U(2^{m+1}):=U(2^{2^{i}})V(2^{2^{i}})+
 V(2^{2^{i}})U(2^{2^{i}})+U(2^{2^{i}})U(2^{2^{i}})+P$.  Certainly,
 conditions 1, 3, 5, 6, 7 hold (and 2 does not apply to this case), and
 condition 5 holds by the observation from the beginning of the proof
 of this theorem.  We claim that condition 4 holds. Indeed, $\widebar
 {f}_{j}\in P\subseteq U(2^{n+1})$ and $g_{j}\subseteq U(2^{n+1})$ for
 every $1\leq j\leq s$, so each $f_{j}\in U(2^{n+1})$. Therefore
 $F_{i}\subseteq U(2^{n+1}) = U(2^{2^{i}})$, as required.

Finally, to finish Case 3, consider the case that $i\notin Z$. In this
 case, take any two monomials $q_{1}, q_{2}$ from
 $V(2^{2^{i}})V(2^{2^{i}})$ and set $V(2^{n+1}):=Kq_{ 1}+Kq_{2}$. Let
 $Q$ be a $K$-linear subspace of $V(2^{2^{i}})V(2^{2^{i}})$ such that
 $Kq_{1}+Kq_{2}+Q=V(2^{2^{i}})V(2^{2^{i}})$ and $Q\cap
 (Kq_{1}+Kq_{2})=0$.  Now set $V(2^{n+1})=Kq_{1}+Kq_{2}$,
 $U(2^{m+1})=U(2^{2^{i}})V(2^{2^{i}})+
 V(2^{2^{i}})U(2^{2^{i}})+U(2^{2^{i}})U(2^{2^{i}})+Q$.  It is easy to
 check all conditions, as in previous cases, noting that
 $2^{n+1}=2^{2{^i}}$ but $i\notin Z$, so that condition 4 holds trivially.
\end{proof} 

\section{The ideal $E$}

The algebra we require will be presented as a factor algebra
$\widebar{A}/E$ for an ideal $E$ that we now define.

\begin{definition}\label{defe} 
{\rm  
 Let $r\in H(n)$ for some $n$, and let $m$ be the natural number such
 that $2^{m}\leq n<2^{m+1}$.  We say that $r\in E(n)$ if and only if
 for all $0\leq j\leq 2^{m+2}-n$ we have 
$$H(j)rH(2^{m+2}-j-n)\subseteq
 U(2^{m+1})H(2^{m+1})+ H(2^{m+1})U(2^{m+1}).$$  We define
 $E=E(1)+E(2)+ \ldots $.}
\end{definition} 

Of course, it is not obvious from this definition that $E$ is an ideal
of $\widebar{A}$. This is the content of the next theorem.

\begin{theorem}
The set $E$ is a two-sided ideal of $\widebar A$.
\end{theorem}

\begin{proof}  
It suffices to show that if $r\in E(n)$ for
some $n$, then $rH(1)\in E(n+1)$ and $H(1)r\subseteq E(n+1)$.
 Let $m$ be the natural number such that 
$2^{m}\leq n<2^{m+1}$.

Consider first the case where  $n<2^{m+1}-1$. Since 
 $r\in E_{n}$ we know that  
$$H(j)rH(2^{m+2}-j-n)\subseteq 
 U(2^{m+1})H(2^{m+1})+H(2^{m+1})U(2^{m+1})$$
 for all $j\leq 2^{m+2}-n$ and this implies that 
$$H(j)rH(1)H(2^{m+2}-j-n-1)\subseteq U(2^{m+1})H(2^{m+1})+
H(2^{m+1})U(2^{m+1})$$ and
 $$H(j)H(1)rH(2^{m+2}-j-n-1)\subseteq U(2^{m+1})H(2^{m+1})+
 H(2^{m+1})U(2^{m+1})$$
  for all $j\leq 2^{m+2}-n-1$. Consequently
$rH(1)\in E_{n+1}$ and $H(1)r\subseteq E_{n+1}$.

Next, consider the case where $n=2^{m+1}-1$. We have to show that
$$H(j)H(1)rH(2^{m+3}-j-n-1)\subseteq U(2^{m+2})H(2^{m+2})+
 H(2^{m+2})U(2^{m+2})$$ 
and  that 
$$H(j)rH(1)H(2^{m+3}-j-n-1)\subseteq
 U(2^{m+2})H(2^{m+2})+H(2^{m+2})U(2^{m+2})$$
for all $j\leq
 2^{m+3}-n-1$.  Consequently we have to show that
$$H(j)rH(2^{m+3}-n-j)\subseteq
U(2^{m+2})H(2^{m+2})+H(2^{m+2})U(2^{m+2})$$
   for all $j\leq 2^{m+3}-n$.

There are three possibilities to consider:
\begin{enumerate}
\item $j+n\leq 2^{m+2}$.
\item  $j\geq 2^{m+2}$.
\item $2^{m+2}<j+n<2^{m+2}+n$.
\end{enumerate}

\Q{\bf Case 1}.  Suppose that $j+n\leq 2^{m+2}$.  
Since $r\in E(n)$ we have
$$H(j)rH(2^{m+2}-n-j)\in U(2^{m+1})H(2^{m+1})
+H(2^{m+1})U(2^{m+1})\subseteq U(2^{m+2}),$$ 
by Theorem~\ref{7props}(6). 
Consequently $H(j)rH(2^{m+3}-n-j)\in U(2^{m+2})H(2^{m+2})$, as
required.

\Q{\bf Case 2}.  Suppose that  
$j\geq 2^{m+2}$. Since $r\in E(n)$ we have
$$H(j-2^{m+2})rH(2^{m+2}-n-j)\in U(2^{m+1})H(2^{m+1})+
H(2^{m+1})U(2^{m+1})\subseteq U(2^{m+2}),$$ 
by Theorem~\ref{7props}(6).
Consequently $H(j)rH(2^{m+3}-n-j)\in H(2^{m+2})U(2^{m+2})$, as
required.

\Q{\bf Case 3}. Suppose that $2^{m+2}<j+n<2^{m+2}+n$. Since $n=2^{m+1}-1$, 
we obtain  $2^{m+1}+1<j$ and $2^{m+3}-j-n>2^{m+1}$. Observe that 
$$H(j)rH(2^{m+3}-j-n)\subseteq H(2^{m+1})H(t)rH(t')H(2^{m+1})$$
for some $t,t'$, where $t+t'=2^{m+2}-n$ (recall that $n=\deg r$).

Since $r\in E(n)$, we obtain 
$$H(t)rH(t')\subseteq U(2^{m+1})H(2^{m+1})
 +H(2^{m+1})U(2^{m+1}).$$ 
Consequently, 
\begin{eqnarray*}H(j)rH(2^{m+3}-j-n)&=&
 H(2^{m+1})H(t)rH(t')H(2^{m+1})\\
&=&
 H(2^{m+1})[U(2^{m+1})H(2^{m+1})+H(2^{m+1})U(2^{m+1})]H(2^{m+1})
 \\
&\subseteq &
U(2^{m+2})H(2^{m+2})+H(2^{m+2})U(2^{m+2}),\end{eqnarray*} 
because 
$H(2^{m+1})U(2^{m+1})\subseteq U(2^{m+2})$ and
 $U(2^{m+1})H(2^{m+1})\subseteq U(2^{m+2})$, by Theorem~\ref{7props}(6).
\end{proof}

\section{
 Definition of  $S(j), W(j), R(j), Q(j)$}

In Section 2, the sets $U(*)$ and $V(*)$ were only defined at powers
of $2$. In this section we define corresponding sets at all other
natural numbers $j$. These are defined in terms of the $U(2^n)$ and
$V(2^n)$ for terms occuring in the binary expansion of $j$. 

Let $j$ be a natural number . Write $j$ in binary form as
 $$j=2^{p_{0}}+2^{p_{1}}+\ldots +2^{p_{n}}$$
 with $0\leq p_{0}<p_{1}<\ldots <p_{n}$. 

\Q
Define $$W(j):=V(2^{p_{0}})V(2^{p_{1}})\ldots 
V(2^{p_{n}})=\prod_{i=0}^{n}V(2^{p_{i}})$$
and set 
 $$S(j):=\sum_{k=0}^{n}S(j,k),$$ with 
 $$S(j,0):=U(2^{p_{0}})H(j-2^{p_{0}})\quad{\rm  and} \quad  
S(j,k):=H(t_{k})U(2^{p_{k}})H(m_{k})$$ where 
 $$t_{k}=\sum_{i=0}^{k-1}2^{p_{i}} \quad{\rm  and} \quad 
m_{k}=\sum_{i=k+1}^{n}2^{p_{i}}$$
for each $j,k>0$.

\Q 
In a similar way, define 
$$Q(j):=V(2^{p_{n}})V(2^{p_{n-1}})\ldots V(2^{p_{0}})
=\prod_{i=0}^{n}V(2^{p_{n-i}})$$
and set  $$R(j):=\sum_{k=0}^{n}R(j,k),$$ with 
 $$R(j,0):=H(j-2^{p_{0}})U(2^{p_{0}})\quad{\rm  and} \quad 
R(j,k)=H(m_{k})U(2^{p_{k}})H(t_{k})$$ where
   $$t_{k}=\sum_{i=0}^{k-1}2^{p_{i}}\quad{\rm  and} \quad 
m_{k}=\sum_{i=k+1}^{n}2^{p_{i}}$$
  for each $j,k>0$.

Note that $W(2^n) = Q(2^n) = V(2^n)$ and that $S(2^n) = R(2^n) = U(2^n)$. 

\begin{lemma}\label{sdirectw}
Let $j$ be a natural number. Then 
\begin{enumerate} 
\item 
$S(j)+W(j)=H(j)$  and  $S(j)\cap
W(j)=0$

\item 
$R(j)+Q(j)=H(j)$  and $R(j)\cap
Q(j)=0$.
\end{enumerate} 
\end{lemma}

\begin{proof} Note that $S(j)\subseteq H(j)$ and $W(j)\subseteq H(j)$
for all $j$. Since $V(2^{p_{i}})+U(2^{p_{i}} )=H(2^{p_{i}})$ for all
 $i$ by Theorem~\ref{7props}(5), we get $S(j)+W(j)=H(j)$.  Observe that
 $V(2^{p_{i}})\cap U(2^{p_{i}} )=0$ for all $i$ by Theorem~\ref{7props}(5).
 Therefore $S(j)\cap W(j)=0$.

The proof of the second claim is similar.
\end{proof} 

\begin{lemma} \label{r=rh+hr}
Let $j$ be a natural number, and let 
$j=2^{p_{0}}+2^{p_{1}}+\ldots +2^{p_{n}}$ be the binary form of $j$
with $0\leq p_{0}<p_{1}<p_{2}<\ldots <p_{n}$.  Let $0<t<n$ and let
 $m=2^{p_{t}}+2^{p_{t+1}}+\ldots +2^{p_{n}}$ and $m'=
 2^{p_{0}}+2^{p_{1}}+\ldots +2^{p_{t-1}}$. Then
 $R(j)=R(m)H(m')+H(m)R(m')$.
 \end{lemma}

\begin{proof} Notice that $m'+m=j$.
  Let $R(j)=\sum_{i=0}^{n}R(j,k)$ be as in the definition above.
 Then $$R(j,k)=H(m_{k})U(2^{p_{k}})H(l_{k})$$ where
   $$l_{k}=\sum_{i=0}^{k-1}2^{p_{i}} \quad{\rm and}\quad 
m_{k}=\sum_{i=k+1}^{n}2^{p_{i}}.$$

Suppose that  $k<t$; so that  $m_{k}\geq m$. Then 
 $R(j,k)=H(m)H(m_{k}-m)U(2^{p_{k}})H(l_{k}).$
  Observe now that $m'=\sum_{i=0}^{t-1}2^{p_i}$
  is the binary form 
  of $m'=j-m$. Therefore $R(m',k)=H(m_{k}-m)U(2^{p_{k}})H(l_{k})$
  for $k<t$. Hence $R(j,k)=H(m)R(m',k)$ for $k<t$, and 
  consequently $$\sum_{i=0}^{t-1}R(j,k)=H(m)R(m').$$

Now suppose that  $k\geq t$; so that  $l_{k}\geq m'$.  Then 
 $$R(j,k)=H(m_{k})U(2^{p_{k}})H(l_{k}-m')H(m'),$$ 
and, arguing as above, 
 $R(j,k)=R(m,k-t+1)H(m')$.  Therefore,
 $$\sum_{i=t}^{n}R(j,k)=R(m)H(m').$$ The result follows.
\end{proof} 

\begin{theorem}\label{rjht}
 For all natural numbers $j,t$  we have
$$
R(j)H(t)\subseteq R(j+t) \quad {\rm and}\quad  H(t)S(j)\subseteq S(j+t).
$$
 \end{theorem}

\begin{proof} 
It is sufficient to show that for every $j$ we have $R(j)H(1)\subseteq
 R(j+1)$ and $H(1)S(j)\subseteq S(j+1)$.  We will show that that
 $R(j)H(1)\subseteq R(j+1)$. The proof that $H(1)S(j)\subseteq S(j+1)$
 is similar.

First, consider the case where 
$j=2^{p+1}-1$ for some $p\geq 0$. Then
 $j=2^{0}+2^{1}+2^{2}+\ldots +2^{p}$.  Consequently
 $R(j)=\sum_{k=0}^{p}R(j,k)$, where
  $$R(j,k)=H(2^{p+1}-2^{k+1})U(2^{k})H(2^{k}-1).$$
   Notice that $R(j+1)=R(2^{p+1})=U(2^{p+1})$. Therefore, it 
   suffices to show that  $R(j,k)H(1)\subseteq R(j+1)=U(2^{p+1})$,
   for every $k\geq 0$.
   Notice that,  $$R(j,k)H(1)=H(2^{p+1}-2^{k+1})U(2^{k})
   H(2^{k}-1)H(1)=H(2^{p+1}-2^{k+1})U(2^{k})H(2^{k}).$$
    Hence, $$R(j,k)H(1)\subseteq
   H(2^{p+1}-2^{k+1})U(2^{k+1}),
$$
by Theorem~\ref{7props}(6).
   Since $H(2^{t})U(2^{t})\subseteq U(2^{t+1})$, 
again by Theorem~\ref{7props}(6), 
we obtain 
   $$H(2^{p+1}-2^{t})U(2^{t})\subseteq H(2^{p+1}-2^{t+1})U(2^{t+1}).$$
   Applying this observation several times
 for $t=k+1$, $t=k+2, \ldots , t=p,$ we get that
 $R(j,k)H(1)\subseteq U(2^{p+1})$, as required.

Next, assume that $j\neq 2^{p+1}-1$ for all  $p$. 
Write $j$  in binary form: $j=2^{p_{0}}+2^{p_{1}}+\ldots +
 2^{p_{n}}$ for some $0\leq p_{0}<p_{1}<p_{2}<\ldots +<p_{n}$.

  First, assume that $p_{0}\neq 0$.  Then
   $j+1=2^{0}+2^{p_{0}}+2^{p_{1}}+\ldots +2^{p_{n}}$ is the binary
   form of $j+1$.  Let $R(j)=\sum_{i=0}^{n}R(j,k)$ and
   $R(j+1)=\sum_{i=0}^{n+1}R(j+1,k)$ be as in the definition.  Now we
   see that $R(j,k)H(1)\subseteq R(j+1,k+1)$.  Therefore
   $R(j)H(1)\subseteq R(j+1)$ as required.

   Next, assume that $p_{0}=0$, and let $t$ be minimal such that
$p_{t}-p_{t-1}>1$.  Then $p_{i}=i$ for all $0\leq i \leq t-1$ and
$p_{t}>t$. Therefore $j=2^{t}-1 +\sum_{i=t}^{n}2^{p_{i}}$.  
By using 
Lemma~\ref{r=rh+hr}, observe that $R(j)=R(m)H(m')+H(m)R(m')$ where
$m'=\sum_{i=0}^{t-1}2^{p_{i}} = 2^t -1$ and $m=\sum_{i=t}^{n}2^{p_{i}}$.
Thus, 
$$
R(j)H(1)=R(m)H(m')H(1)+
H(m)R(m')H(1).
$$ 
Since $m'=2^{t}-1$, we get $R(m')H(1)\subseteq
R(m'+1)=R(2^{t})$, by the first part of the proof. Therefore
 $$
R(j)H(1)\subseteq R(m)H(m'+1)+H(m)R(m'+1).
$$
 
Observe that the binary form of $j+1$ is  
$j+1=2^{t}+ \sum_{i=t}^{n}2^{p^{i}}$.
 Recall that $m=\sum_{i=t}^{n}2^{p^{i}}$ and that $2^{t}=m'+1$.
   Now from Lemma~\ref{r=rh+hr}, 
we get $$R(j+1)=R(m)H(m'+1)+H(m)R(m'+1).$$ 
   Consequently $R(j)H(1)\subseteq R(j+1)$, and the lemma follows.
\end{proof} 

\section{Estimation of the Gelfand-Kirillov dimension}

In order to estimate the Gelfand-Kirillov dimension of
$\widebar{A}/E$, we need to recognise when certain homogenous elements
are in $E$. The next theorem provides a sufficient condition for this
to happen.

\begin{theorem}\label{sh+hr} 
 If $r\in H(n)$ for some $n$ and $r\in S(t)H(n-t)+H(t)R(n-t)$
 for all $0\leq t\leq n$, then $r\in E$.
\end{theorem}

\begin{proof} 
Suppose that $r\in S(t)H(n-t)+H(t)R(n-t)$
 for all $0\leq t\leq n$. 
Let $m$ be the natural number such that 
$2^{m}\leq n<2^{m+1}$. By the definition of $E$ we 
have to show that  for all
$j\leq 2^{m+2}-n$ we have
$$H(j)rH(2^{m+2}-j-n)\subseteq U(2^{m+1})H(2^{m+1})+
H(2^{m+1})U(2^{m+1}).$$
 Consider the three possibilities:
 \begin{itemize}
 \item[1.] $j+n>2^{m+1}$ and $2^{m+1}-j\geq 0$.
 \item[2.] $j+n>2^{m+1}$ and $2^{m+1}-j<0$.
 \item[3.]$j+n\leq 2^{m+1}$.
 \end{itemize}

  \Q{\bf Case 1.} Suppose first that  $j+n>2^{m+1}$ and $2^{m+1}-j\geq 0$. 
  Set $t=2^{m+1}-j$. By assumption, 
 $r\in  S(t)H(n-t)+H(t)R(n-t)$, since $0\leq t$ and $n-t=n+j-2^{m+1}>0$.
 Therefore 
 $$
H(j)rH(2^{m+2}-j-n)\subseteq H(j)
\left\{
S(t)H(n-t)  +H(t)R(n-t)
\right\}
H(2^{m+2}-j-n).
$$
 Now, since $j+t=2^{m+1}$, we get
 $$H(j)rH(2^{m+2}-j-n)\subseteq H(j)S(t)H(2^{m+1})+H(2^{m+1})
 R(n-t)H(2^{m+1}-n+t).$$ 
By Theorem~\ref{rjht}, $H(j)S(t)\subseteq
 S(2^{m+1})=U(2^{m+1})$. Similarly, $R(n-t)H(2^{m+1}-n+t)\subseteq
 R(2^{m+1})=U(2^{m+1})$, by Theorem~\ref{rjht}. Therefore
 $$
H(j)rH(2^{m+2}-j-n)\subseteq  U(2^{m+1})H(2^{m+1})+H(2^{m+1})U(2^{m+1}),
$$ 
as required.

  \Q{\bf Case 2.} Suppose that $j+n>2^{m+1}$ and $j>2^{m+1}$. Then
$j=2^{m+1}+b$ for some $b> 0$.  Since $j+n\leq 2^{m+2}$, it follows that 
$b+n\leq 2^{m+1}$; and $b\leq 2^{m}$, 
since $n\geq 2^{m}$.
Now take $t
=2^{m}-b$. Then $0\leq t\leq n$.  Hence, by assumption, 
$r\in S(t)H(n-t)+H(t)R(n-t)$. Consequently,
$$
H(j)rH(2^{m+2}-j-n)\subseteq H(j)\left\{
S(t)H(n-t)+H(t)R(n-t)
\right\}
H(2^{m+2}-j-n).
$$
Since $j=2^{m+1}+b$ we obtain 
$$
H(j)rH(2^{m+2}-j-n)
\subseteq H(2^{m+1})H(b)\left\{
S(t)H(n-t)+H(t)R(n-t)
\right\} 
H(2^{m+2}-j-n).
$$

Consider the two terms that occur on the right hand side of this
containment separately. 

First, consider the term
$
H(2^{m+1})H(b)
S(t)H(n-t)
H(2^{m+2}-j-n)$

Note that  $t+b=2^{m}$ and that $j+t = 2^{m+1} + 2^m$. 
Hence,  $H(b)S(t)\subseteq S(2^{m}) = U(2^m)$, 
by Theorem~\ref{rjht}; and so  
 \begin{eqnarray*}
H(2^{m+1})H(b)
S(t)H(n-t)
H(2^{m+2}-j-n) &=& 
H(2^{m+1})U(2^m)H(2^{m+2}-j-t)\\ 
&=&H(2^{m+1})U(2^m)H(2^m) 
\\&\subseteq& H(2^{m+1})U(2^{m+1}),
\end{eqnarray*} 
as required. 

Next, consider the term 
$
H(2^{m+1})H(b)
H(t)R(n-t)
H(2^{m+2}-j-n)
$.

Observe that $R(n-t)H(2^{m+2}-j-n)\subseteq R(2^{m}) = U(2^m)$, by
Theorem~\ref{rjht}, since $j+t=2^{m+1}+2^{m}$. Also, $H(b)H(t) = H(t+b)
= H(2^m)$.

Hence, 
\begin{eqnarray*}
H(2^{m+1})H(b)
H(t)R(n-t)
H(2^{m+2}-j-n) 
   &\subseteq& 
H(2^{m+1})H(2^m)U(2^m) \\
&\subseteq& H(2^{m+1})U(2^{m+1}),
\end{eqnarray*}
as required.

  Consequently,  
  $H(j)rH(2^{m+2}-j-n)\subseteq H(2^{m+1})U(2^{m+1})$, as required.

  \Q{\bf Case 3.}
Suppose that  $j+n\leq 2^{m+1}$. Then $j\leq 2^{m}$, since $n\geq 2^{m}$.
 Set $t:=2^{m}-j$. Then $0\leq t\leq n$. By assumption,
 $r\in S(t)H(n-t)+H(t)R(n-t)$. Note that
 $j+t=2^{m}$. Therefore, 
$$H(j)rH(2^{m+2}-n-j)\subseteq H(j)\left\{ 
S(t)H(n-t)+H(t)R(n-t)
\right\} 
H(2^{m+2}-n-j).$$
Theorem~\ref{rjht} gives $H(j)S(t) \subseteq S(j+t) = S(2^m) = U(2^m)$; 
so that 
$$H(j)rH(2^{m+2}-n-j)\subseteq U(2^{m})H(2^{m+2}-j-t)+
 H(2^{m})R(n-t)H(2^{m+2}-n-j).$$ 
  
Consider the two terms on the right hand side of this containment. 

First, 
$$
U(2^{m})H(2^{m+2}-j-t)= 
U(2^m)H(2^{m+2} - 2^m) 
= U(2^m)H(2^m)H(2^{m+1})= U(2^{m+1})H(2^{m+1}). 
$$

Secondly, 
note that $2^m - n+t = 2^m -n + (2^m -j) =2^{m+1} - (n+j) \geq 0$ and that 
$2^{m+2} - n-j = 2^{m+1} + 2^{m+1} -n-j = 2^{m+1} + (2^m -n+t)$; and so 
 \begin{eqnarray*}
H(2^{m})R(n-t)H(2^{m+2}-n-j)
&=& 
H(2^{m})R(n-t)H(2^{m}-n+t)H(2^{m+1})\\ 
&\subseteq& H(2^m)R(2^m)H(2^{m+1}) \\
&=& H(2^m)U(2^m)H(2^{m+1})\subseteq U(2^{m+1})H(2^{m+1}),
\end{eqnarray*}
as required. 

Consequently $H(j)rH(2^{m+2}-n-j)\subseteq U(2^{m+1})H(2^{m+1})$.
 This finishes the proof.
\end{proof}

%%%%%%%%%%%%%%%%%%%%%%%%%%%%%%%%%%%%%%%%%%%%%%%%%%%%%%

We can now estimate the size of the subspaces $Q(n)$ and $W(n)$. 

\begin{theorem}\label{dimqn}
For all $n>0$, we have $\dim_{K}Q(n)\leq 3^{17}n^{9}$ and $\dim_{K}W(n)\leq
3^{17}n^{9}$.
\end{theorem}

\begin{proof}
We will show that $\dim_{K}Q(n)\leq 3^{17}n^{9}$. The proof
that $\dim_{K}W(n) \leq 3^{17}n^{9}$ is similar. Let   
$n=2^{p_{0}}+2^{p_{1}}+\ldots +2^{p_{m}}$ with $0\leq
p_{0}<p_{1}<p_{2}<\ldots <p_{m}$. 
Then
$Q(n)=V(2^{p_{0}})V(2^{p_{1}})\ldots V(2^{p_{m}})$.  
Consequently
$\dim_{K}Q(n)=\prod_{i=0}^{m}\dim_{K}V(2^{p_{i}})$.  
Therefore
$\dim_{K}Q(n)\leq \prod_{i=0}^{\lfloor\log (n)\rfloor}\dim_{K}V(2^{i})$, 
where
$\lfloor\log(n)\rfloor$ 
is the largest integer not exceeding $\log(n)$. 
Recall that, from  Theorem 2, either $\dim_{K}V(2^{i})=2$ 
or $i\in S$, where $S=\{[2^{i}-i-1,
2^{i}-1] \mid  i=5,6, \ldots \}$.  

Now let $c_{i}=\prod_{t=2^{i}-i-1}^{2^{i}-1}\dim_{K}V(2^{t})$.  We see
that $c_{i}=\prod_{j=0}^{i}\dim_{K}V(2^{2^{i}-i-1+j})
=\prod_{j=0}^{i}2^{2^{j}}<2^{2^{i+1}}$, by Theorem ~\ref{7props}(2).  Since
$\dim_{K}Q(n)\leq \prod_{i=0}^{\lfloor \log
(n)\rfloor}\dim_{K}V(2^{i})$, we have $\dim_{K}Q(n)\leq cc'$ where
$c=\prod_{i\notin S,i\leq \lfloor\log (n)\rfloor}\dim_{K}V(2^{i})$ and
$c'= \prod_{i\in S, i\leq \lfloor\log (n)\rfloor }\dim_{K}V(2^{i})$.
First, observe that $c \leq 2^{\lfloor\log (n)\rfloor+1}\leq 2n$.
Next, let $q$ be the maximal number such that $2^{q}-q-1\leq
\lfloor\log (n)\rfloor$.  Then $c'\leq \prod_{i=0}^{q}c_{i}\leq
\prod_{i=0}^{q}2^{2^{q+1}}\leq 2^{2^{q+2}}$.  Observe that
$2^{q-1}\leq 2^{q}-(q+1)$, if $q\geq 3$.  Therefore
$2^{2^{q+2}}=2^{2^{q-1}2^3}=(2^{2^{q-1}})^{8}\leq \left(
2^{\lfloor\log (n)\rfloor}\right)^8 \leq n^8$, provided that $q\geq 3$.
Observe that if $n\geq 16$ then the maximal number $q$ with
$2^{q}-q-1\leq \lfloor\log (n)\rfloor$ is indeed greater than or equal
to $3$; so that $\dim_K(Q(n) \leq 2n.n^8 = 2n^9$ provided that $n\geq
16$.  If $n<16$, then $\dim_K(Q(n)) \leq \dim_K(H(n)) < 3^{16}$, since
$Q(n) \subseteq H(n)$ for each $n$ and $A$ is generated by the three
elements $x, y, z$.  
Therefore $\dim_{k}Q(n)\leq 3^{17}n^{9}$ for each $n$.
\end{proof} 
 
%%%%%%%%%%%%%%%%%%%%%%%%%%%%%%%%%%%%%%%%%%%%%%%%%%%%%%%

After all this preparation, we can now estimate the Gelfand-Kirillov
dimension of our factor algebra.

\begin{theorem}
$\gkdim(\widebar {A}/E)\leq 20$.
\end{theorem}

\begin{proof}
Let $0\leq j\leq n$. 
Note that $S(j)+W(j)=H(j)$ and $R(n-j)+Q(n-j)=H(n-j)$. 
It follows that 
$$
H(n) = H(j)H(n-j) = W(j)Q(n-j) + \{ S(j)H(n-j) + H(j)R(n-j)\}.
$$
Thus, 
$$
\dim \left(\frac{H(n)}{\{ S(j)H(n-j) + H(j)R(n-j)\}}\right) 
\leq \dim W(j)Q(n-j) \leq 
(3^{17}n^9)^2 = 3^{34}n^{18}.
$$
Let 
$$
\theta:H(n) \longrightarrow 
\bigoplus_{j=0}^{n}\, \frac{H(n)}{\{ S(j)H(n-j) + H(j)R(n-j)\}}
$$
be the natural map. 
Then 
$$
\ker(\theta) 
= \{r\in H(n) \mid r\in S(j)H(n-j) + H(j)R(n-j) 
{\rm ~for~each~} 0\leq j\leq n \} \subseteq E(n)
$$
by Theorem~\ref{sh+hr}. 
Thus, 
$$\dim\left(\frac{H(n)}{E(n)}\right) 
\leq \dim\left(\frac{H(n)}{\ker(\theta)} \right) 
\leq 3^{34}n^{18}(n+1).
$$
Consequently, $\gkdim(\widebar{A}/E) \leq 20$.

\end{proof}

%%%%%%%%%%%%%%%%%%%%%%%%%%%%%%%%%%%%%%%%%

\section{$\widebar{A}/E$ is nil but not nilpotent}

It remains to show that the algebra $\widebar{A}/E$ is nil but not
nilpotent.  We show that $\widebar{A}/E$ is nil by showing that the
elements ${f_{i}}^{10w_{i}}$ defined in Section~\ref{one} belong to
$E$. In order to see that $\widebar{A}/E$ is not nilpotent we show
that the $K$-subspaces $V(2^n)$ are not contained in $E$.

%%%%%%%%%%%%%%%%%%%%%%%%%%%%%%%%%%%%%%%%%%

\begin{lemma}\label{bwifi}
Let $Z$, $\{f_{i}\}_{i\in Z}$, $\{F_{i}\}_{i\in Z}$ be as in 
Theorem~\ref{existFi}. 
 Fix any $i\in Z$ and suppose that $m+2>2^{i}$. Then $B_{w_{i}}(F_{i})\cap
H(2^{m+2})\subseteq U(2^{m+1})H(2^{m+1})+ H(2^{m+1})U(2^{m+1})$.
\end{lemma}

\begin{proof} 
We know that $F_{i}\subseteq H(2^{2^{i}})$ and 
$w_{i}=4r_{i}$ where $r_{i}=2^{2^{i}}$. Set $w:=w_{i}$ and $r:=r_{i}$.
By the assumptions of this lemma
$2^{m+1}\geq 2^{2^{i}}$. 
Observe that 
$B_{w}(F_{i})\subseteq B_{r}(F_{i})$, since $w=4r$. 
Also,  $B_{r}(F_{i})\subseteq B_{r}(U(r))$, 
since $F_{i}\subseteq U(2^{2^{i}})=U(r)$ by Theorem~\ref{7props}(4). 
  Consequently $$B_{w}(F_{i})\subseteq B_{r}(U(r)).$$
 Therefore, it is sufficient to show that
$$
B_{r}(U_{r})\cap H(2^{m+2})\subseteq 
U(2^{m+1})H(2^{m+1})+H(2^{m+1})U(2^{m+1})
$$
for all $m$ such that $2^{m+1}\geq 2^{2^{i}}=r$ . 
 We will proceed  by induction on $m$.
If $m+1=2^{i}$ then $2^{m+2}=2r$;  so that 
 $B_{r}(U(r))\cap H(2r)=U(r)H(r)+H(r)U(r)$, by the definition 
 of $B_{r}(U(r))$, and the fact that $U(r)\subseteq H(r)$.

Suppose now that the result holds for some $m$, 
with $2^{m+1}\geq
2^{2^{i}}=r$. We will prove that the result holds 
for $m+1$.  We have to show that
$$
B_{r}(U(r))\cap H(2^{m+3})\subseteq U(2^{m+2})H(2^{m+2})+
H(2^{m+2})U(2^{m+2}).
$$ 
Observe that, since $r$ divides $2^{m+2}$, we obtain 
$$
B_{r}(U(r))\cap H(2^{m+3})=\left\{B_{r}(U(r))\cap
H(2^{m+2})\right\}H(2^{m+2}) 
+H(2^{m+2})\left\{B_{r}(U(r))\cap H(2^{m+2})\right\},
$$ 
by the definition of $B_{r}(U(r))$.  By the induction assumption
$$
B_{r}(U(r))\cap H(2^{m+2})\subseteq
U(2^{m+1})H(2^{m+1})+H(2^{m+1})U(2^{m+1})\subseteq U(2^{m+2})
$$ 
by
Theorem~\ref{7props}(6). Hence,  $B_{r}(U(r))\cap H(2^{m+3})\subseteq
U(2^{m+2})H(2^{m+2})+H(2^{m+2})U(2^{m+2})$ and the result
follows.\end{proof}

%%%%%%%%%%%%%%%%%%%%%%%%%%%%%%%%%%%%%%%%%%%%%%%%%%%

\begin{theorem}
Let $Z$, $\{f_{i}\}_{i\in Z}$ be as in Lemma 1. Let $i\in Z$ and
let $I$ be the two sided ideal of $R$ generated by $f_{i}^{10w_{i}}$
where $w_{i}=4.2^{2^{i}}$. Then $I\subseteq E$.
\end{theorem}

\begin{proof} 
Let $r\in I$. Then $r=\sum_{p=10w_{i}}^{s}r_{p}$ for some $r_{p}\in
H(p)$, and some $s$.  Fix $n$, with $10w_{i}\leq n\leq s$.  
It is sufficient to 
show that $r_{n}\in E$. Let $m$ be the natural number such that 
$2^{m}\leq n<2^{m+1}$.  Note that $10w_i = 40.2^{2^{i}}$; so 
$40.2^{2^{i}}\leq n<2^{m+1}$.  Hence
$m+1>2^{i}$. In order to  show that $r_{n}\in E$, we have to show that 
$$H(j)r_{n}H(2^{m+2}-n-j)\subseteq
U(2^{m+1})H(2^{m+1}) 
+H(2^{m+1})U(2^{m+1}),$$ 
for
every $0\leq j\leq 2^{m+2}-n$. 

Now $r\in I$ yields
 $H(j)rH(2^{m+2}-n-j)\in I$. Consequently,
 $H(j)rH(2^{m+2}-n-j)\subseteq B_{w_{i}}(F_{i})$, by Theorem~\ref{existFi}; 
so
 $$H(j)\left(\sum_{p=10w_{i}}^{s}r_{p}\right)H(2^{m+2}-n-j)
 \subseteq B_{w_{i}}(F_{i}).$$
 
It follows that 
$H(j)r_{p}H(2^{m+2}-n-j)\subseteq B_{w_{i}}(F_{i})$,
 for every $p$ with $10w_i \leq p\leq s$, since   $B_{w_{i}}(F_{i})$ is 
homogeneous and $r_{p}\in H(p)$ for  every $p$.

In particular, 
 $$H(j)r_{n}H(2^{m+2}-n-j)\subseteq B_{w_{i}}(F_{i}).$$ 
Now,  since
 $H(j)r_{n}H(2^{m+2}-n-j)\subseteq H(2^{m+2})$ and $m+2>2^{i}$, we
 have 
\begin{eqnarray*}
H(j)r_{n}H(2^{m+2}-n-j) &\subseteq& 
B_{w_{i}}(F_i)\cap H(2^{m+2})\\ &\subseteq& H(2^{m+1})U(2^{m+1})
 +U(2^{m+1})H(2^{m+1}),
\end{eqnarray*} 
 by Lemma~\ref{bwifi}, and this completes the proof.  
\end{proof} 

%%%%%%%%%%%%%%%%%%%%%%%%%%%%%%%%%%%%%%%%%%%%%

 The next two results are now immediate. 

\begin{corollary}
 Let $Z$, $\{f_{i}\}_{i\in Z}$ be as in Lemma 1. Let $N$ be the
 two-sided ideal in $A$ generated by elements from the set
 $\{f_{i}^{10w_{i}}\}_{i\in Z}$ where $w_{i}=4.2^{2^{i}}$.  Then
 $N\subseteq E$.
\end{corollary}

%%%%%%%%%%%%%%%%%%%%%%%%%%%%%%%%%%%%%%%%%%%%%%%%%

\begin{theorem}
The algebra $\widebar{A}/E$ is a nil algebra.
\end{theorem}
\begin{proof} 
This follows  from the previous theorem  and Lemma 1. 
\end{proof} 

%%%%%%%%%%%%%%%%%%%%%%%%%%%%%%%%%%%%%%%%%%%%%%%%
Finally, we show that $\widebar{A}/E$ is not nilpotent. 

\begin{theorem}
The algebra $ \widebar{A}/E$ is not nilpotent.
\end{theorem}

\begin{proof} 
Recall that $V(2^{n+1})\subseteq V(2^{n})V(2^{n})$, for every $n>0$, by
Theorem~\ref{7props}(7). It follows easily, by induction, that $V(2^m)
\subseteq V(2)^{2^{m-1}}$. Thus, it is sufficient to show that $V(2^m)
\not\subseteq E$.  

Recall that, by Theorem~\ref{7props}(3), $V(2^{m})$ is generated by
monomials, for all $m$.  Therefore, there are $0\neq r\in V(2^{m})$ and 
$0\neq r'\in H(2^{m+2}-2^{m})$ such that $0\neq rr'\in V(2^{m+2})$.  Suppose
that $r\in E$; so that, in fact, $r\in E(2^m)$. By using the defining
property of $E$, see Definition~\ref{defe}, with $j=0$ and $n=2^m$, we obtain
\begin{eqnarray*} 
0\neq rr' &\in & H(0)rH(2^{m+2} -0-2^m) 
\subseteq U(2^{m+1})H(2^{m+1}) + H(2^{m+1})U(2^{m+1})
\subseteq U(2^{m+2}).
\end{eqnarray*} 
Thus, $0 \neq rr'\in V(2^{m+2}) \cap U(2^{m+2}) = 0$, 
contradicting Theorem~\ref{7props}(5). 

Hence, $r\not\in E$; so that $V(2^m) \not\subseteq E$, as required. 
\end{proof} 

%%%%%%%%%%%%%%%%%%%%%%%%%%%%%%%%%%%
In conclusion, we have proved: 

\begin{theorem}
The finitely generated algebra $ \widebar{A}/E$  is nil, but not nilpotent, 
and has Gelfand-Kirillov dimension not exceeding $20$. 
\end{theorem}

%%%%%%%%%%%%%%%%%%%%%%%%%%%%%%%%%%%%%%%%%%%%%%%%%%%

\end{document}